\documentclass[11pt]{amsart}
\usepackage{latexsym, amsthm, amsmath, bbm}
\usepackage{amssymb}
\usepackage{amsfonts}
\usepackage{amstext}
\usepackage{multicol}

\newtheorem{theorem}{Theorem}[section]
\newtheorem{corollary}[theorem]{Corollary}
\newtheorem{lemma}[theorem]{Lemma}
\newtheorem{proposition}[theorem]{Proposition}
\newtheorem{conjecture}[theorem]{Conjecture}
\newtheorem{definition}[theorem]{Definition}

\theoremstyle{remark}
\newtheorem{remark}[theorem]{Remark}

\theoremstyle{remark}

\theoremstyle{remark}
\newtheorem{example}[theorem]{Example}

\numberwithin{equation}{section}

\newcommand{\ep}{\epsilon}
\newcommand{\del}{\delta}

\newcommand{\Z}{{\mathbb{Z}}}

\newcommand{\E}{{\mathbb{E}}}
\newcommand{\reals}{\mathbb{R}}
\newcommand{\ints}{\mathbb{Z}}
\newcommand{\nats}{\mathbb{N}}

\newcommand{\rats}{\mathbb{Q}}

\newcommand{\ess}{\mathcal{S}}

\newcommand{\Zmn}{\mathbb{Z}_m^n}

\newcommand{\inv}{^{-1}}

\newcommand{\ovl}{\overline}

\newcommand{\til}{\widetilde}

\newcommand{\into}{\hookrightarrow}

\newcommand{\subeq}{\subseteq}

\newcommand{\bslash}{\backslash} 

\newcommand{\mbf}{\mathbf}
\def\diam{\operatorname{diam}}

\def\dist{\operatorname{dist}}

\def\XXint#1#2#3{{\setbox0=\hbox{$#1{#2#3}{\int}$} 
\vcenter{\hbox{$#2#3$}}\kern-.5\wd0}}

\begin{document}
\title{Spaces of small metric cotype}

\author{E. Veomett}
\address{Department of Mathematics and Computer Science, California State University, East Bay,  Hayward, California 94542-3092}
\email{ellen.veomett@csueastbay.edu}

\author{K. Wildrick}
\address{Department of Mathematics and Statistics, University of Jyv\"askyl\"a, PL 35 MaD, 40014 Jyv\"askyl\"an yliopisto, Finland}
\thanks{The second author was supported by the Academy of Finland grants 120972 and 128144}
\email{kevin.wildrick@jyu.fi}

\date{}

\begin{abstract}
Mendel and Naor's definition of metric cotype extends the notion of the Rademacher cotype of a Banach space to all metric spaces. Every Banach space has metric cotype at least $2$. We show that any metric space that is bi-Lipschitz equivalent to an ultrametric space has infimal metric cotype $1$.  We discuss the invariance of metric cotype inequalities under snowflaking mappings and Gromov-Hausdorff limits, and use these facts to establish a partial converse of the main result. 
\end{abstract}
\maketitle
\section{Introduction}
\noindent

Recent connections between theoretical computer science, geometric functional analysis, and analysis on metric spaces have led to significant progress in these fields \cite{CKN}, \cite{Markov}, \cite{MN}.  Type and cotype inequalities, which  play an important role in the geometry of Banach spaces, provide an obstruction to embedding metric spaces into highly structured Banach spaces where efficient algorigthms may be available \cite{TypeCoAndK}, \cite{Kalton}. In \cite{MN}, Mendel and Naor introduced a satisfactory notion of cotype for general metric spaces, and used it to solve a variety of problems regarding geometric embeddings. Among other important results, Mendel and Naor establish a non-linear version of the Maurey-Pisier theorem for spaces with no metric cotype.  However, several results in \cite{MN} apply only to Banach spaces, though the statements are sensible in a general metric setting.  The prospect of extending these results is tantalizing. In this paper, we begin this process by examining situations in which the theory of metric cotype differs from its linear counter-part.

Recall that a Banach space $X$ is said to have Rademacher type $p\in [1,2]$ if there exists a number $T\geq 1$ such that for each positive integer $n$ and each $x_1, x_2, \dotsc, x_n \in X$,
\begin{equation*}
\sum_{i=1}^n ||x_i||_X^p \geq T^{-p} \E_\ep||\sum_{i=1}^n \epsilon_i x_i ||_X^p,\end{equation*}
where the expectation is taken over uniformly distributed $\ep \in \{-1,1\}^n.$  
Similarly, $X$ is said to have Rademacher cotype $q \in [2,\infty)$ if there exists a number $C\geq 1$ such that for each positive integer $n$ and each $x_1, x_2, \dotsc, x_n \in X$,
\begin{equation*}
 \sum_{i=1}^n ||x_i||_X^q \leq C^q \E_\ep||\sum_{i=1}^n \epsilon_i x_i ||_X^q \end{equation*}
where the expectation is as before.  

We now describe the notion of metric cotype established in \cite{MN}.  Let $X$ be a metric space, $n\in \nats$, and $m \in 2\nats$. We denote the $n$-fold product of the integers modulo $m$ by $\Zmn$.  For a function $f \colon \Zmn \to X$, an index $1\leq j \leq n$, and a vector $\del \in \{-1,0,1\}^n$, define
$$f_j(\ep)=f\left(\ep+\frac{m}{2}\mathbf{e}_j\right) \ \text{and}\ f_{\del}(\ep)=f(\ep + \del),$$
where $\mathbf{e}_j$ denotes the $j$th standard basis vector.

\begin{definition}[Mendel-Naor]
Let $1\leq p\leq q<\infty$. A metric space $(X,d)$ supports a $(p,q)$-metric cotype inequality if there is a constant $\Gamma \geq 1$ and a scaling function $m_{p,q}\colon \ints \to 2\ints$ with the following property.  Given $n \in \ints $ and setting $m=m_{p,q}(n)$, every mapping $f \colon \Zmn \to X$ satisfies  
\begin{equation}\label{metric cotype} \E_{\ep \in \Zmn} \sum_{j=1}^n d(f(\ep),f_j(\ep))^p \leq \Gamma^pm^pn^{1-(p/q)}\E_{\ep \in \Zmn} \E_{\del \in \{-1,0,1\}^n}d(f(\ep),f_\del(\ep))^p.\end{equation} 
The expectations are taken over uniformly distributed $\ep \in \Z_m^n$ and $\del \in \{-1,0,1\}^n$. If $(X,d)$ supports a $(q,q)$-metric cotype inequality, then we say that it supports a metric cotype $q$ inequality.\end{definition}

Mendel and Naor's main result \cite[Theorem 1.4]{MN} shows the aptness of the above definition.

\begin{theorem}[Mendel-Naor]\label{MN main} If a Banach space supports a $(p,q)$-metric cotype inequality, then it has Rademacher cotype $q'$ for every $q'>q$.  Conversely, a Banach space with Rademacher cotype $q$ supports a $(p,q)$-metric cotype inequality for every $1 \leq p\leq q$.
\end{theorem}

Motivated by these results, given a metric space $(X,d)$, we denote
$$q_X = \inf\{q \geq 1 : \text{there is $1 \leq p \leq q$ such that $X$ supports a $(p,q)$-metric cotype inequality}\}.$$
In this notation, Theorem \ref{MN main} implies that if $X$ is a Banach space, then $q_X$ is the infimum over $q$ such that $X$ has Rademacher cotype $q$.  

Recall that if a Banach space has Rademacher cotype $q$, then it has Rademacher cotype $q'$ for all $q' \geq q$, and that no Banach space supports a Rademacher cotype inequality with exponent less than $2$.  Our main result shows that certain highly disconnected metric spaces support ``better" metric cotype inequalities than are possible for Banach spaces. 

A metric space $(X,d)$ is an \emph{ultrametric space} if for all triples $x,y,z \in X$, 
$$d(x,y) \leq \max\{d(x,z),d(z,y)\}.$$
Ultrametric spaces are fundamental mathematical objects, used in theoretical computer science \cite{Ramsey}, $p$-adic analysis \cite{UMCalc}, and mathematical biology \cite{Evolution}.  Informally, they are spaces that arise as the leaves of a metric tree. 

Metric cotype inequalities are naturaly invariant under bi-Lipschitz mappings; see Proposition \ref{scaled bi-Lip invariance}. Hence, in this paper, we consider metric spaces that are bi-Lipschitz equivalent to an ultrametric space. Examples of such spaces include all finite metric spaces, the standard middle third Cantor set, and the space $\{2^{-i}\}_{i \in \ints} \subeq \reals$. Every separable ultrametric space isometrically embeds into the Hilbert space of square summable sequences \cite{Timan}.  Hence, every metric space that is bi-Lipschitz equivalent to a separable ultrametric space supports a metric cotype $2$ inequality.  However, more is true.

\begin{theorem}\label{STSTheorem} Suppose that a metric space $(X,d)$ is bi-Lipschitz equivalent to an ultrametric space.  Then $(X,d)$ supports a metric cotype $q$ inequality for every $q>1$.  In particular, $q_X =1$.  
\end{theorem}

The main tool in the proof of Theorem \ref{STSTheorem} is an edge isoperimetric inequality on a graph that corresponds to the right hand side of the metric cotype inequality \eqref{metric cotype}.  This graph can be thought of as an $\ell^{\infty}$-discrete torus.  We derive the needed inequalities from similar results on the $\ell^1$-discrete torus in \cite{MR1137765}.

The second portion of this paper discusses the invariance of metric cotype under various types of mappings between general metric spaces, and under limiting processes.  We use the results regarding the invariance of cotype to work towards a converse of Theorem \ref{STSTheorem}, employing the formalism of $s$-snowflake spaces established by Tyson and Wu \cite{TW}.  Let $s \in [1,\infty]$.  A metric space $(X,d)$ is said to be an \emph{$L^s$-metric space} if for all triples $x,y,z \in X$,
$$d(x,y) \leq \begin{cases} 
			\left(d^s(x,z) + d^s(z,y)\right)^{1/s} & 1\leq s <\infty, \\
			\max\{d(x,z),d(z,y)\} & s=\infty.\\
\end{cases}$$
Note that a metric space is an $L^\infty$-metric space if and only if it is an ultrametric space.   We define the \emph{snowflake index} $s_X$ of a metric space $(X,d)$ by 
$$s_X = \sup\{s : \text{$(X,d)$ is bi-Lipschitz equivalent to an $L^s$-metric space}\}.$$
A metric space $(X,d)$ is said to be an \emph{$s$-snowflake} if $(X,d)$ is bi-Lipschitz equivalent to an $L^{s}$-metric space and $s=s_X$.  Hence, a metric space is an $\infty$-snowflake if and only if it is bi-Lipschitz equivalent to an ultrametric. In this notation, Theorem \ref{STSTheorem} implies that if $(X,d)$ is an $\infty$-snowflake, then $q_X=1$.  

Using results in \cite{TW} due to Laakso, we show the following statements, which are in the direction of a converse to Theorem \ref{STSTheorem}. 

\begin{theorem}\label{1 snowflake} Suppose that $(X,d)$ is a $1$-snowflake.  Then $q_X \geq 2$.  
\end{theorem}

\begin{theorem}\label{p snowflake} Suppose that $(X,d)$ is a $s$-snowflake for some $s\in [1,\infty)$.  If $1 \leq p\leq q<2$, then $(X,d)$ does not support a $(p,q)$-metric cotype inequality in which $m_{p,q}(n) =O(n^{1/q})$. 
\end{theorem}

As mentioned in \cite{MN} and \cite{GMN}, it is of interest to know if the conclusion of Theorem \ref{p snowflake} implies that $q_X \geq 2$.  Also, it would be interesting to know if the hypothesis in Theorem \ref{p snowflake} can be weakened to the condition that $s_X < \infty$, and if the hypothesis in Theorem \ref{STSTheorem} can be weakened to the condition $s_X = \infty$. 

Comparing Theorems \ref{STSTheorem} and \ref{p snowflake} leads us to the following conjecture. Similar phenomena have been observed for conformal dimension; see \cite{LeoConfDim}.

\begin{conjecture}\label{q gap} There is no metric space $(X,d)$ for which $1<q_X<2$.  
\end{conjecture}

It seems likely that a similar analysis could be made of the corresponding notion of metric type given in 
\cite{EnfloType}.

We would like to thank Leonid Kovalev for inspiration, Assaf Naor for encouragement and helpful suggestions, and Jeremy Tyson for pointing out the result of Laakso, Theorem \ref{line fitting 1 snowflake}. Also, thanks to Daniel Meyer, Jaime Radcliffe, and the excellent referees for critical comments. 

\section{Classes of mappings}

We now establish notation for the various classes of mappings we will consider throughout this paper. Let $\phi \colon (Y,d_Y) \to (X,d_X)$ be a mapping between metric spaces. 

For $c>0$ and $L \geq 1$, we say that $\phi$ is an \emph{$c$-scaled $L$-bi-Lipschitz embedding} if for all $x,y \in Y$,
$$\frac{d_Y(x,y)}{L}\leq cd_X(\phi(x),\phi(y)) \leq Ld_Y(x,y).$$
Such a mapping distorts absolute distances by a fixed multiplicative factor, adjusted for scaling. If $c=1$, such a map is simply called an \emph{$L$-bi-Lipschitz embedding}.

For $c>0$, $\alpha >0$, and $L \geq 1$, we say that $\phi$ is an \emph{$c$-scaled $(\alpha,L)$-snowflaking embedding} if for all $x,y \in Y$
$$\frac{d_Y(x,y)^\alpha}{L}\leq cd_X(\phi(x),\phi(y)) \leq Ld_Y(x,y)^\alpha.$$
If $c=1$, such a map is simply called an \emph{$(\alpha, L)$-snowflaking embedding}. When $0<\alpha\leq 1$, such a mapping may be thought of as a scaled bi-Lipschitz mapping defined on the snowflaked space $(Y,d_Y^\alpha)$.  The term ``snowflake" comes from the fact that the standard parameterization of the Von Koch snowflake curve by the unit circle is a $(\log_43,L)$-snowflaking embedding for some $L\geq 1$.  

We will often supress the value of the scaling factor $c$ in the above classes of mappings, as it usually plays no quantitative role. 

Let $\eta\colon [0,\infty) \to [0,\infty)$ be a homeomorphism. We say that $\phi$ is an \emph{$\eta$-quasisymmetric} embedding If for each triple of distinct points $x,y,z \in Y$,
$$\frac{d_X(\phi(x),\phi(y))}{d_X(\phi(x),\phi(z))} \leq \eta\left(\frac{d_Y(x,y)}{d_Y(x,z)}\right).$$
Such a mapping distorts relative distances by a controlled amount.  Note that every snowflaking embedding is quasisymmetric. Quasisymmetric embeddings are a generalization of conformal mappings to the metric space setting.  For an introduction to their basic properties, see \cite{LAMS}.   

\section{$\infty$-snowflakes}\label{um section}
\noindent

We first provide a characterization of $\infty$-snowflakes similar to that described in \cite[Section 2]{TW} and \cite[Section 14.24]{LAMS}.  Throughout, we denote the cardinality of a set $A$ by $|A|$. 

\begin{proposition}\label{tree metric}
Let $(X,d)$ be a metric space. Then the following conditions are equivalent:
\begin{itemize}
\item[(i)]the space $(X,d)$ is $L$-bi-Lipschitz equivalent to an ultrametric space,
\item[(ii)] the space $(X,d)$ is $\eta$-quasisymmetrically equivalent to an ultrametric space,
\item[(iii)] there is a constant $C \geq 1$ such that for every subset $S$ of $X$ satisfying $1 < |S| < \infty$, there is a subset $\emptyset \neq A \subsetneq S$ such that 
$$\dist(A, S\bslash A) \geq \frac{\diam S}{C}.$$
\end{itemize}
\end{proposition}

\begin{remark}\label{quantitative} The above proposition is quantitative in the following sense.  If condition (i) is satisfied, then condition (ii) is satisfied with $\eta(t)=L^2t.$  If condition (ii) is satisfied, then condition (iii) is satisfied with $C=2\eta(1)$. If condition (iii) is satisfied, then condition (i) is satisfied with $L=C$.  
\end{remark}

A space satisfying condition (iii) above is said to have the \emph{$C$-finite separation property}.    We prove Proposition \ref{tree metric} via three easy lemmas. We leave the proof of the first to the reader, and the proof of the third is essentially the same as the proof of \cite[Proposition~15.7]{Broken}, so we omit it as well.

\begin{lemma}\label{FSP qs invariant} If $f\colon X \to Y$ is an $\eta$-quasisymmetric homeomorphism, and $X$ has the $C$-finite separation property for some $C\geq 1$, then $Y$ has the $2\eta(C)$-finite separation property. 
\end{lemma}


\begin{lemma}\label{um has fsp} If $(X,d)$ is an ultrametric space, then it has the $1$-finite separation property.
\end{lemma}

\begin{proof}  Let $S$ be a subset of $X$ satisfying $1<|S| <\infty$.  Set $D=\diam S$, and let $\mathcal{D}$ be a maximal collection of points in $S$ such that $d(x,y) \geq D$ for all $x,y\in \mathcal{D}$. Since $S$ is a finite set, the definition of diameter implies that $\mathcal{D}$ has at least two points. Fix any $x_0\in \mathcal{D},$ and define 
$$A = \{x \in S: d(x_0,x)<D\}.$$
Let $x\in A$ and $y \in S\bslash A$.  Since $\mathcal{D}$ is maximal, there exists $y_0 \in \mathcal{D}$ such that $d(y_0,y)<D$. Noting that
$$D\leq d(x_0,y_0)\leq \max\{d(x_0,x),d(x,y_0)\},$$
we see that $d(x,y_0)\geq D.$  Similarly, 
$$D\leq d(x,y_0) \leq \max\{d(x,y),d(y,y_0)\},$$
and so $d(x,y)\geq D$.  This implies the desired result.
\end{proof}

\begin{lemma}\label{FSP implies biLip to um} Suppose that a metric space $(X,d)$ satisfies the $C$-finite separation property for some $C \geq 1$.  Then there is an ultrametric space that is $C$-bi-Lipschitz equivalent to $X$.
\end{lemma}

%

\begin{proof}[Proof of Proposition \ref{tree metric}] An $L$-bi-Lipschitz mapping is $\eta$-quasisymmetric with $\eta(t)=L^2t$.  Thus (i) implies (ii).  Lemmas \ref{FSP qs invariant} and \ref{um has fsp} show that (ii) implies (iii).  Lemma \ref{FSP implies biLip to um} shows that (iii) implies (i).  
\end{proof}

%

The conditions in Proposition \ref{tree metric} provide a natural tree structure on any bounded $\infty$-snowflake. Let $\ess$ denote the collection of finite sequences $\alpha =(\alpha_i)_{i=1}^k$ such that $\alpha_i \in \{0,1\}$ for each index $i$. We also consider the empty sequence $\emptyset$ to be in $\ess$. Given $\alpha \in \ess$ and $\del \in \{0,1\}$, we set
$$\alpha\frown\del=(\alpha_1,\hdots,\alpha_k,\del) \in \ess \ \text{and}\ \del\frown \alpha = (\del,\alpha_1,\hdots,\alpha_k) \in \ess. $$ 
 We say that $\alpha =(\alpha_i)_{i=1}^k \in \ess$ is an \emph{ancestor} of $\alpha'=(\alpha'_i)_{i=1}^{k'} \in \ess$ if $k\leq k'$ and $\alpha_j=\alpha'_j$ for $j\leq k$.  We also call $\alpha'$ a \emph{descendant} of $\alpha$. 

\begin{definition}\label{CSTS}
A $C$-separated tree structure, $C \geq 1$, on a metric space $(X,d)$ is a collection $\{A_\alpha\}_{\alpha \in \ess}$ of subsets of $X$ such that 
\begin{enumerate}
\item \label{empty} $A_\emptyset = X,$
\item \label{separate} given distinct points $x$ and $y$ in $X$, there are sequences $\alpha$ and $\beta$ in $\ess$ such that neither is an ancestor of the other, and such that $x \in A_\alpha$ and $y \in A_\beta$,  
\end{enumerate}
and for any $\alpha \in \ess$,
\begin{enumerate}
\setcounter{enumi}{2}
\item \label{partition} $A_{\alpha}=A_{\alpha \frown 0}\cup A_{\alpha \frown 1}$,
\item \label{non-trivial} If $|A_{\alpha}| > 1$, then $A_{\alpha \frown 0}\neq \emptyset \neq A_{\alpha \frown 1},$
\item \label{gap} If $|A_{\alpha}|> 1$, then $\diam A_\alpha \leq C \dist(A_{\alpha\frown 0}, A_{\alpha \frown 1}).$
\end{enumerate}
\end{definition}

\begin{lemma}\label{finite equiv} Let $C\geq 1$ and let $(X,d)$ be a finite metric space.  Then $(X,d)$ supports a $C$-separated tree structure if and only if it has the $C$-finite separation property. 
\end{lemma}

\begin{proof} 
Suppose that $\{A_\alpha\}_{\alpha \in \ess}$ is a $C$-separated tree structure on $(X,d)$, and let $S$ be a finite subset of $X$ with at least two points.  By properties \eqref{empty} and \eqref{separate} of Definition \ref{CSTS}, there is a sequence $\beta \in \ess$ with the property that $S$ is contained in $A_\beta$ but not in any of its descendants.  Then 
\begin{align*}\dist(S \cap A_{\beta \frown 0}, S\bslash A_{\beta \frown 0}) &\geq \dist(A_{\beta \frown 0}, A_{\beta \frown 1})\\ &\geq \frac{\diam A_\beta}{C} \\ & \geq \frac{\diam S}{C}, \end{align*}
and so $(X,d)$ has the $C$-finite separation property.

We now show that if $(X,d)$ has the $C$-finite separation property, then it supports a $C$-separated tree structure. The assertion is trivially true when $|X|\leq 1$.  So suppose that $|X| > 1$, and that the assertion is true for any metric space $Y$ such that $|Y| <|X|$.  Define $A_\emptyset =X$.  By the $C$-finite separation property, we may find subset $S \subeq X$ such that $0 < |S| <|X|$ such that 
$$\dist(S,X\bslash S) \geq \frac{\diam X}{C}.$$
Set $A_{0}=S$ and $A_{1}=X\bslash S$. Since the $C$-finite separation property is inherited by subsets, the induction hypothesis provides $C$-separated tree structures $\{A'_\alpha\}_{\alpha \in \ess}$ and $\{A''_\alpha\}_{\alpha \in \ess}$ on $A_0$ and $A_1$ respectively.  Given a non-empty sequence $\alpha \in \ess$, we may write $\alpha = \del \frown \beta$ for some $\beta \in \ess$ and $\del \in \{0,1\}$.  If $\del = 0$, define $A_\alpha = A'_\beta$, and if $\del = 1$, define $A_\alpha = A''_\beta$. Then $\{A_\alpha\}_{\alpha \in \ess}$ is a $C$-separated tree structure on $X$.  
\end{proof}

We may now characterize $\infty$-snowflakes in terms of separated tree structures.

\begin{proposition}\label{um tree char} Let $(X,d)$ be a metric space.  Then there is a number $L \geq 1$ such that $(X,d)$ is $L$-bi-Lipschitz equivalent to an ultrametric space if and only if there is a number $C \geq 1$ such that every finite subspace of $(X,d)$ supports a $C$-separated tree structure.  Moreover, $L$ and $C$ depend only on each other.
\end{proposition}
\begin{proof} Suppose that $(X,d)$ is $L$-bi-Lipschitz equivalent to an ultrametric space, $L \geq 1$.  By Proposition \ref{tree metric}, the space $(X,d)$ has the $C$-finite separation property for some $C \geq 1$ depending only on $L$.  Hence every finite subspace of $(X,d)$ also has the $C$-finite separation property.  Lemma \ref{finite equiv} now shows that every finite subspace of $(X,d)$ supports a $C$-separated tree structure.

If every finite subspace of $(X,d)$ supports a $C$-separated tree structure, then $(X,d)$ has the $C$-finite separation property.  Proposition \ref{tree metric} now shows that $(X,d)$ is $L$-bi-Lipschitz equivalent to an ultrametric space for some $L \geq 1$ depending only on $C$.
\end{proof}

We will use the following basic result regarding separated tree structures, the simple proof of which we leave to the reader.

\begin{lemma}\label{earlier branches}
Let $\{A_{\alpha}\}_{\alpha \in \ess}$ be a $C$-separated tree structure.  If $\alpha$ and $\alpha'$ are elements of $\ess$ such that neither is an ancestor of the other, and neither of $A_\alpha$ and $A_{\alpha'}$ are empty, then  
\begin{equation*}
\dist(A_{\alpha}, A_{\alpha'}) \geq \frac{\max\{\diam A_\alpha, \diam A_{\alpha'}\}}{C}.\end{equation*}
\end{lemma}


\section{The proof of Theorem \ref{STSTheorem}}\label{main section}

For a positive integer $n$ and an even positive integer $m$, we denote by $L^n_m$ the graph with vertex set $\Z_m^n$ and edge set 
$$E_L = \{(\ep,\ep'): \ep' = \ep+\frac{m}{2}\mathbf{e}_j \text{ for some } j=1, 2, \dotsc, n\}.$$
Similarly, we denote by $R^n_m$ the graph with vertex set $\Z_m^n$ and edge set
$$E_R = \{(\ep,\ep'): \ep' =\ep +\del \text{ for some } 0 \neq \delta \in \{-1, 0, 1\}^n\}.$$
  Given a set $A \subeq \Zmn$, we define the \emph{edge boundary of $A$ in $E_R$}  by
\begin{equation}\label{boundary def}\partial_R{A} = \{(\ep, \ep') \in E_R:  |\{\ep, \ep'\} \cap A| = 1 \}.\end{equation}

The following edge isoperimetric inequality for $R_m^n$ allows us to give a lower bound for the right-hand side of \eqref{metric cotype}.  Note that it is of Euclidean type.  We derive it from a similar statement on the $\ell_1$-discrete torus \cite{MR1137765}.

\begin{theorem}\label{boundary}
Suppose that $A \subset \Z_m^n$ satisfies $|A|\leq \frac{m^n}{2}$.  Then 
$$|\partial_R{A}| \geq 2|A|^{(n-1)/n}.$$
\end{theorem}

\begin{proof}
Consider the graph $T=(\Zmn, E)$, where a pair of vertices $(x_1, x_2, \dots, x_n) \in \Z_m^n$ and $(y_1, y_2, \dots, y_n) \in \Z_m^n$ form an edge in $E$ if there exists an index $i \in \{1,\hdots, n\}$ for which 
\begin{equation*}
|x_j-y_j| = \begin{cases} 1 & j=i, \\
			0 & j \neq i.\\
	\end{cases}\end{equation*}
Suppose $A \subeq \Zmn$ satisfies $|A| \leq \frac{m^n}{2}$, and define
$$\partial_TA = \{(x,y) \in E: |\{x,y\} \cap A | = 1\}.$$
In \cite[Theorem 8]{MR1137765},  Bollob\'{a}s and Leader show that 
\begin{equation}\label{BLTheorem}
|\partial_TA| \geq \min\{2|A|^{1-1/r}rm^{n/r-1} : r = 1, 2, \dots, n\}
\end{equation}
We note that $T$ and $R_m^n$ have the same vertex set, and $E \subeq E_R$. Thus,
\begin{equation*}
|\partial_RA| \geq | \partial_TA|.
\end{equation*}
We may write $|A| = c^n$ for some $c \in \reals$ satisfying $c < m$.  By \eqref{BLTheorem}, 
\begin{align*}
|\partial_TA| &\geq \min\{2c^{n(1-1/r)}rm^{n/r-1} : r = 1, 2, \dots, n\} \\
& \geq  \min\{2c^{n(1-1/r)}m^{n/r-1} : r = 1, 2, \dots, n\} \\
& \geq 2c^{n-1} = 2|A|^{(n-1)/n},
\end{align*}
as desired. \end{proof}

\begin{proof}[Proof of Theorem \ref{STSTheorem}]  We assume that $(X,d)$ is $L$-bi-Lipschitz equivalent to an ultrametric space, $L \geq 1$.  Let $q>1$, and let $m_{q,q} \colon \ints \to 2\ints$ be any function such that 
$$m_{q,q}(n)  \geq (n3^n)^{1/(q-1)}.$$
Fix $n \in \ints$ and set $m=m_{q,q}(n)$.  We will show that for any $f \colon \Zmn \to X$, 
\begin{equation}\label{metric cotype q} \E_{\ep \in \Zmn} \sum_{j=1}^n d(f(\ep),f_j(\ep))^q \leq \Gamma^qm^q \E_{\ep \in \Zmn} \E_{\del \in \{-1,0,1\}^n}d(f(\ep),f_\del(\ep))^q,\end{equation} 
for some $\Gamma \geq 1$ that depends only on $L$. 

By Proposition \ref{um tree char}, there is a $C$-separated tree structure $\{A_{\alpha}\}_{\alpha \in \ess}$ on $f(\Zmn)$, where $C \geq 1$ depends only on $L$.  We may find $i_0 \in \nats$ such that if the length of $\alpha$ is greater than or equal to $i_0$, then $A_{\alpha}$ is either empty or a singleton. For each $\alpha \in \ess$, set $F_{\alpha} = f\inv(A_\alpha)$, with the convention that $F_{\alpha}=\emptyset$ if $A_{\alpha}=\emptyset$.  

The graph $L^n_m$ is related to the left-hand side of the metric cotype inequality \eqref{metric cotype q} in the following way: 
\begin{equation}\label{LHS interpret}
\E_{\ep \in \Zmn} \sum_{j=1}^n d(f(\ep),f_j(\ep))^q = 2m^{-n}\sum_{(\ep,\ep') \in E_L} d(f(\ep),f(\ep'))^q.
\end{equation}
The factor of $2$ occurs because the graph $L^n_m$ is not directed.

We estimate the quantity in \eqref{LHS interpret} by partitioning $E_L$ into sets where uniform estimates are possible.  For an integer $i \geq 0$, we denote by $\mbf{1}_i \in \ess$ the string consisting of $i$ ones, and set 
$$E_L^i=\{(\ep,\ep') \in E_L: \ \{\ep,\ep'\} \subeq F_{\mbf{1}_i} \ \text{and}\ \{\ep,\ep'\}\cap F_{\mbf{1}_i \frown 0} \neq \emptyset\}.$$
Since there are $n$ edges attached to any vertex in $L^n_m$, we estimate that $|E_L^i|\leq n|F_{\mbf{1}_i\frown 0}|$, and so 
$$\sum_{(\ep,\ep') \in E^i_L} d(f(\ep),f(\ep'))^q \leq  n|F_{\mbf{1}_i\frown 0}|(\diam A_{\mbf{1}_i})^q.$$
Moreover, $\{E_L^i\}_{i=0}^{i_0}$ covers $E_L$, and $\diam A_{\mbf{1}_{i_0}} = 0.$ Thus, \eqref{LHS interpret} implies that 
\begin{equation}\label{LHS est}\E_{\ep \in \Zmn} \sum_{j=1}^n d(f(\ep),f_j(\ep))^q \leq m^{-n} \sum_{i=0}^{i_0-1} 2n|F_{\mbf{1}_i\frown 0}|(\diam A_{\mbf{1}_i})^q.\end{equation}

We now give a lower bound for the right-hand side of \eqref{metric cotype q}. We relate it to the graph $R^n_m$ via the identity
\begin{equation}\label{RHS interpret} 
\E_{\ep \in \Zmn} \E_{\del \in \{-1,0,1\}^n}d(f(\ep),f_\del(\ep))^q = 2(3m)^{-n}\sum_{(\ep,\ep') \in E_R} d(f(\ep),f(\ep'))^q.
\end{equation}
An edge in $E_R$ appears in at most two of the sets $\{\partial_R F_{\mbf{1}_i\frown 0}\}_{i=0}^{i_0-1}.$ Hence 
\begin{equation}\label{RHS decomp} \sum_{(\ep,\ep') \in E_R} d_X(f(\ep),f(\ep'))^q \geq \frac{1}{2} \sum_{i=0}^{i_0-1} \sum_{(\ep,\ep') \in \partial_R F_{\mbf{1}_i \frown 0}} d_X(f(\ep),f(\ep'))^q. \end{equation}

Let $(\ep,\ep') \in \partial_R F_{\mbf{1}_i\frown 0}$ for some integer $0\leq i < i_0$.  We claim that 
\begin{equation}\label{RHS help} d_X(f(\ep),f(\ep')) \geq \frac{\diam A_{\mbf{1}_i}}{C}.\end{equation}
Assume without loss of generality that $\ep \in F_{\mbf{1}_i \frown 0}$ and $\ep' \notin  F_{\mbf{1}_i\frown 0}$.  If $\ep' \in F_{\mbf{1}_{i+1}}$, then \eqref{RHS help} follows from the definition of a $C$-separated tree structure.  If $\ep'\notin F_{\mbf{1}_{i+1}}$, then $\ep' \notin F_{\mbf{1}_i}$. This implies that we may find $\alpha \in \ess$ such that $\ep' \in F_{\alpha}$ and neither $\alpha$ nor $\mbf{1}_i$ is an ancestor of the other. Since $\ep \in F_{\mbf{1}_i}$, Lemma \ref{earlier branches} now yields \eqref{RHS help}.  

The inequalities \eqref{RHS interpret}, \eqref{RHS decomp}, and \eqref{RHS help} show that 
\begin{equation}\label{RHS est} C^q m^q \E_{\ep \in \Zmn} \E_{\del \in \{-1,0,1\}^n}d(f(\ep),f_\del(\ep))^q \geq m^{-n}\sum_{i=0}^{i_0-1} 3^{-n}m^q|\partial_R F_{\mbf{1}_i\frown 0}|\left(\diam A_{\mbf{1}_i}\right)^q.\end{equation}

The desired inequality \eqref{metric cotype q} with $\Gamma=C$ will now follow from \eqref{LHS est} and \eqref{RHS est} provided we additionally show that given $0\leq i <i_0$,
\begin{equation}\label{calculation}3^{-n}m^q|\partial_R F_{\mbf{1}_i, 0}| \geq 2n|F_{\mbf{1}_i, 0}|.\end{equation}

For each integer $0 \leq i < i_0$, the sets $F_{\mbf{1}_i\frown 0}$ and $F_{\mbf{1}_{i+1}}$ are disjoint subsets of $\Zmn$.  Thus, by relabeling the separated tree structure if needed, we may assume without loss of generality that for each integer $0 \leq i < i_0$, the set $F_{\mbf{1}_i\frown 0}$ has cardinality no greater than $m^n/2.$   Hence by Theorem \ref{boundary}, for each integer $0 \leq i < i_0$,
$$|\partial_R F_{\mbf{1}_i\frown 0}| \geq 2|F_{\mbf{1}_i\frown 0}|^{(n-1)/n}.$$ 
The definition of $m_{q,q}$ also implies that $3^{-n}m^q \geq mn$.  Hence, using the above isoperimetric inequality and the trivial estimate $|F_{\mbf{1}_i\frown 0}| \leq m^n$, we see that
$$3^{-n}m^q |\partial_R F_{\mbf{1}_i\frown 0}| \geq \left(2n|F_{\mbf{1}_i\frown 0}|\right)\left(m|F_{\mbf{1}_i\frown 0}|^{-1/n}\right) \geq  2n|F_{\mbf{1}_i\frown 0}|,$$
as desired. \end{proof}

We record the following corollary, which one might consider to be a very weak non-linear version of the Maurey-Pisier theorem for cotype $1$. In light of Conjecture \ref{gap}, we suspect that a stronger statement is true.

\begin{definition} Let $\ep>0$, and let $a,b \in X$.  An \textit{$\ep$-chain} in $X$ connecting $a$ to $b$ is a finite sequence of points $x_0=a, x_1,\hdots,x_k=b$ such that for each $i=0,\hdots,k-1$,
$$d(x_i,x_{i+1})<\ep.$$
\end{definition}

\begin{corollary}\label{very weak MP} Suppose that $(X,d)$ is a metric space with $q_X>1$.  Then for all $C \geq 1$, there are distinct points $a, b \in X$ and a $d(a,b)/C$-chain connecting $a$ to $b$ with diameter $d(a,b)$.  
\end{corollary} 

\begin{proof}  Fix $C \geq 1.$  Let $\mathcal{F}$ be the family of finite subsets $S \subeq X$ such that if $\emptyset \neq A \subsetneq S$, then 
\begin{equation}\label{no sep}\dist(A,S\bslash A) < \frac{\diam S}{C}.\end{equation}
Then by Theorem \ref{STSTheorem} and Proposition \ref{tree metric}, the family $\mathcal{F}$ is non-empty.

Let $S \in \mathcal{F}$ and set $M:= \diam S$.  Since $S$ is finite, we may find points $a,b \in S$ such that $d(a,b)=M$.  Consider the set
$$\mathcal{C}:= \{x \in S: \ \text{there exists an $(M/C)$-chain in $S$ connecting $a$ to $x$}\}.$$
Clearly $\mathcal{C}$ is non-empty as it contains $a$. If $\mathcal{C}\neq S$, then \eqref{no sep} provides points $x \in \mathcal{C}$ and $y \in S\bslash \mathcal{C}$ such that $d(x,y)<M/C$.  Adding the point $y$ to the $(M/C)$-chain connecting $a$ to $x$ yields an $(M/C)$-chain connecting $a$ to $y$, a contradiction.  Thus $\mathcal{C}=S$, and so there is an $(M/C)$-chain in $S$ connecting $a$ to $b$.  
\end{proof}

\section{The invariance of metric cotype}

Metric cotype inequalities are invariant under scaled bi-Lipschitz embeddings.  We leave the proof of the following statement to the reader.

\begin{proposition}\label{scaled bi-Lip invariance} Let $1\leq p\leq q<\infty$ and $L \geq 1$.  Suppose $X$ and $Y$ are metric spaces such that there is a scaled $L$-bi-Lipschitz embedding of $Y$ into $X$.  If $X$ supports a $(p,q)$-metric cotype inequality with constant $\Gamma$ and scaling function $m_{p,q}(n)$, then $Y$ supports a $(p,q)$-metric cotype inequality with the same scaling function $m_{p,q}(n)$ and constant depending only on $L$ and $\Gamma$. 
\end{proposition}

A key feature of cotype in the Banach space setting is that it is preserved by a variety of much larger classes of non-linear mappings than just bi-Lipschitz embeddings.  The following theorem from \cite{Note} is an example of such a result. 

\begin{theorem}[Naor]\label{qs Banach invariance} Suppose that $f \colon W \into V$ is a quasisymmetric embedding of Banach spaces, and assume that $V$ has non-trivial type.  Then 
\begin{align*} \inf\{q: \text{W supports a metric cotype}\ & q\ \text{inequality}\} \\ & \leq \inf \{q: \text{V supports a metric cotype $q$ inequality}\}.\end{align*}
In particular, $q_W \leq q_V$.  
\end{theorem}

The proof of Theorem \ref{qs Banach invariance} relies heavily on the classical Maurey-Pisier theorem \cite{MP} and the equivalence between metric and Rademacher cotype.  The assumption that $V$ have non-trivial type is equivalent to the assumption that $V$ is $K$-convex.  As shown in \cite[Theorem~4.1]{MN}, this implies that in considering metric cotype $q$ inequalities on $V$, one may assume that $m_{p,q}(n) = O(n^{1/q})$.  This is also a crucial fact in the proof.  Theorem \ref{qs Banach invariance} and the equivalence between Rademacher and metric cotype imply that $L_p$ does not quasisymmetrically embed in $L_2$ if $p>2$.  We do not know if Theorem \ref{qs Banach invariance} is valid in a general setting.  This seems to be an interesting and difficult question.

The class of scaled snowflaking embeddings, which is smaller than the class of quasisymmetric embeddings, preserves the existence of \emph{some} $(p,q)$-metric cotype inequality \emph{provided that} $m_{p,q}(n)=O(n^{1/q})$.  

\begin{proposition}\label{sf invariance} Let $1\leq p\leq q<\infty$, and suppose that $\phi \colon Y \into X$ is  a $c$-scaled $(\alpha, L)$-snowflaking embedding.  If $X$ supports a $(p,q)$-metric cotype inequality with constant $\Gamma$ and a scaling function satisfying $m_{p,q}(n)\leq Kn^{1/q}$ for some $K \geq 1$, then $Y$ supports an $(\alpha p,q)$-metric cotype inequality with the same scaling function $m_{p,q}$, and constant depending only on $\alpha$, $L$, $\Gamma$, and $K$.  \end{proposition}

\begin{proof} Set $p' =\alpha p$.  Let $\Gamma \geq 1$ and $m_{p,q}$ be the constant and scaling function associated to the $(p,q)$-metric cotype inequality on $X$, respectively.  Let $n \in \nats$ and set $m=m_{p,q}(n)$, and consider any $f \colon \Zmn \to Y$. Then
\begin{align*}\mathbb{E}_{\ep \in \Zmn}\sum_{j=1}^n  d_Y(f(\ep),f_j(\ep))^{p'} &\leq (Lc)^p\mathbb{E}_{\ep \in \Zmn}\sum_{j=1}^n d_X(\phi \circ f(\ep),\phi\circ f_j(\ep))^{p} \\
&\leq (Lc)^p\Gamma^{p} m^pn^{1-(p/q)}\E_{\ep \in \Zmn} \E_{\del \in \{-1,0,1\}^n} d_X(\phi\circ f_\del(\ep),\phi\circ f(\ep))^p \\
& \leq L^{2p}\Gamma^{p} m^pn^{1-(p/q)}\E_{\ep \in \Zmn} \E_{\del \in \{-1,0,1\}^n}d_Y(f(\ep), f_{\del}(\ep))^{p'}. \end{align*}
It follows from \cite[Lemma 2.3]{MN} that $m \geq \Gamma\inv n^{1/q}$, and by assumption there is a constant $K \geq 1$ such that $m \leq Kn^{1/q}$. Thus 
$$m^pn^{1-(p/q)} \leq K^pn =K^pn^{p'/q}n^{1-(p'/q)} \leq K^p\Gamma^{p'}m^{p'}n^{1-(p'/q)}.$$
Combining these estimates yields the desired result, with constant $L^{2p/p'}\Gamma^{(p+p')/p'}K^{p/p'}$ and scaling function $m_{p,q}$.\end{proof}

\begin{remark}  In Theorem \ref{qs Banach invariance}, it is essentially shown that a $(q,q)$-metric cotype inequality is preserved under quasisymmetric embeddings. However, even though the class of snowflaking embeddings is much smaller, in Proposition \ref{sf invariance} the first exponent may change.  It is not clear if a better result is possible in this general setting.  By using H\"older's inequality and similar estimates for sums, if $p=q$ in Proposition~\ref{sf invariance}, then one can show that $Y$ supports the following inequality:
$$\mathbb{E}_{\ep \in \Zmn}\left(\sum_{j=1}^n  d_Y(f(\ep),f_j(\ep))^{q}\right)^{\alpha} \leq 
L^{2q} \Gamma^qm^q\E_{\ep \in \Zmn} \left(\E_{\del \in \{-1,0,1\}^n}d(f(\ep),f_\del(\ep))^q\right)^{\alpha}.$$ 
It would be interesting to know the consequences of such an inequality. 
\end{remark}

The following corollary now follows from the fact that a Hilbert space supports a metric cotype $2$ inequality with $m_{2,2} =O(n^{1/2})$ \cite[Proposition~3.1]{MN}.

\begin{corollary}\label{embedding application} Suppose that $Y$ is a metric space that $(\alpha, C)$-snowflake embeds into a Hilbert space.  Then $Y$ supports a $(2\alpha, 2)$-metric cotype inequality. In particular, $q_Y \leq 2$. \end{corollary}

\begin{remark}\label{Assouad consequence} A metric space $Y$ is doubling if there is a number $N \in \nats$ such that for every point $y \in Y$ and radius $r>0$, the ball $B(y,r)$ may be covered by at most $N$ balls of radius $r/2$.  Assouad's Theorem \cite[Theorem~12.2]{LAMS} states that if $Y$ is doubling, then for any $0 < \alpha <1$, we may find $L \geq 1$ and $n \in \nats$, depending only on $\alpha$ and the doubling constant of $Y$, such that there is an $(\alpha, L)$-snowflaking embedding of $Y$ into $\reals^n$. Thus Corollary \ref{embedding application} implies that for all $0<\alpha<1$, the space $Y$ supports a $(2\alpha, 2)$-metric cotype inequality.  We do not know if this implies that $Y$ supports a metric cotype $2$ inequality. 
\end{remark}

\section{Gromov-Hausdorff limits}

It is not difficult to show that if a metric space $X$ supports a metric cotype inequality, then the completion $\ovl{X}$ satisfies the same inequality, possibly with a different constant.  Essentially the same proof shows that metric cotype inequalities are inherited by Gromov-Hausdorff limits.  In this section we use this and related facts to give lower bounds for $q_X$ in a variety of circumstances.  For an introduction to Gromov-Hausdorff limits, see \cite{Burago}.  

We begin with some notation. For $c>0$, a (possibly non-continuous) function $\phi\colon Y \to X$ of metric spaces is said to be a \emph{$c$-rough isometry} if for all $y,z \in Y$,
$$|d_Y(y,z) - d_X(\phi(y),\phi(z))|\leq c.$$
A subset $A$ of a metric space $X$ is said to be $c$-dense if for all $x \in X$, there is $a \in A$ such that $d_X(a, x) \leq c$.  If the image of a mapping $\phi \colon Y \to X$ is $c$-dense in $X$, then there is a \emph{rough inverse mapping} $\til{\phi}\colon X \to Y$ with the property that $d_X(\phi\circ \til{\phi}(x), x) \leq c$. 

We say that a sequence of compact metric spaces $\{(X_k,d_k)\}$ Gromov-Hausdorff converges to a metric space $(X,d)$ if for every $c>0$, there is $K \in \nats$ so that if $k \geq K$, then there is a $c$-rough isometry from $X_k$ to $X$ with $c$-dense image.  

There is a version of Gromov-Hausdorff convergence that is more appropriate in the non-compact setting.  A \emph{pointed metric space} is a pair $((X,d),p)$ where $(X,d)$ is a metric space and $p$ is a point in $X$.  A sequence of pointed metric spaces $\{((X_k,d_k),p_k))\}_{k \in \nats}$ \emph{pointed Gromov-Hausdorff converges} to a pointed metric space $((X,d),p)$ if for all $r>0$ and all $c>0$, there is  $K \in \nats$ such that if $k \geq K$, then there is a $c$-rough isometry $\phi \colon X_k \to X$ such that $\phi(X_k)$ is $c$-dense in $B(p,r-c)$, and $\phi(p_k)=p$.  This notion restricts to the standard notion in the case that all spaces involved are compact \cite[Exercise~8.1.2]{Burago}.

\begin{proposition}\label{GH converge}  Let $1 \leq p \leq q < \infty$.  Fix $\Gamma \geq 1$ and a scaling function $m_{p,q} \colon \nats \to 2\nats$.  Suppose that each metric space in the sequence of pointed spaces $\{((X_k,d_k),p_k)\}_{k\in \nats}$ supports a $(p,q)$-metric cotype inequality with constant $\Gamma$ with scaling function $m_{p,q}$.  If $\{((X_k,d_k),p_k)\}_{k \in \nats}$ Gromov-Hausdorff converges to a pointed space $((X,d),p_\infty)$, then $(X,d)$ supports a $(p,q)$-metric cotype inequality with constant $4\Gamma$ and scaling function $m_{p,q}$. 
\end{proposition} 

\begin{proof}  Let $m = m_{p,q}(n)$, let $f \colon \Zmn \to X$ be a function.  Since $\Zmn$ is finite, we may find $r>0$ such that $f(\Zmn) \subeq B_X(p_\infty,r/2)$.  Let $0<c<r/2.$ By assumption, there is $K \in \nats$ such that if $k \geq K$, then there is a $c$-rough isometry $\phi \colon X_k \to X$ such that $\phi(X_k)$ is $c$-dense in $B(p_\infty,r-c)$, and hence in $f(\Zmn)$.  As a result there is a $c$-rough inverse $\til{\phi}$ of $\phi$ defined on $f(\Zmn)$.  Thus, for each $\ep \in \Zmn$ and $j =1,\hdots, n$, 
\begin{align*}&|d_X(f(\ep),f_j(\ep)) - d_{X}(\phi \circ \til{\phi} \circ f(\ep), \phi \circ \til{\phi} \circ f_j(\ep))| \leq 2c,\ \text{and}\ \\
			 &|d_{X}(\phi \circ \til{\phi} \circ f(\ep), \phi \circ \til{\phi} \circ f_j(\ep)) - d_{X_k}(\til{\phi} \circ f(\ep),\til{\phi} \circ f_j(\ep))| \leq c.\end{align*}
These facts, along with the elementary inequality $(a+b)^p \leq 2^p(a^p+b^p)$, where $a,b\geq 0$, imply that 
$$ \E_{\ep \in \Zmn} \sum_{j=1}^n d_X(f(\ep),f_j(\ep))^p  \leq n(6c)^p +2^{p}\E_{\ep \in \Zmn} \sum_{j=1}^n d_{X_k}(\til{\phi} \circ f(\ep), \til{\phi} \circ f_j(\ep))^p.$$
We now apply the $(p,q)$-metric cotype inequality to the map $\til{\phi} \circ f$ and make similar estimates, yielding
\begin{align*} \E_{\ep \in \Zmn} \sum_{j=1}^n d_{X_k}(\til{\phi} \circ f(\ep),  \til{\phi} \circ  f_j&(\ep))^p \leq \Gamma^p m^pn^{1-(p/q)} \E_{\ep \in \Zmn} \E_{\del \in \{-1,0,1\}^n}d_{X_k}(\til{\phi} \circ f(\ep),\til{\phi} \circ f_\del(\ep))^p \\
&\leq \Gamma^p m^pn^{1-(p/q)} \left((6c)^p + 2^{p}\E_{\ep \in \Zmn} \E_{\del \in \{-1,0,1\}^n}d_{X}(f(\ep),f_\del(\ep))^p \right). \end{align*}
We now see that 
\begin{align*}\E_{\ep \in \Zmn} \sum_{j=1}^n d_X(f(\ep),f_j(\ep))^p \leq (6c)^p(n&+ 2^{p}\Gamma^p m^pn^{1-(p/q)})\\&+4^{p}\Gamma^p m^pn^{1-(p/q)}  \E_{\ep \in \Zmn} \E_{\del \in \{-1,0,1\}^n}d_{X}(f(\ep),f_\del(\ep))^p.\end{align*}
Letting $c$ tend to zero now yields the desired result.  
\end{proof}

\begin{remark}\label{completion} Let $x_0 \in X$. The definitions imply that the constant sequence $\{((X,d),x_0)\}_{k=1}^\infty$ pointed Gromov-Hausdorff converges to the completed space $(\ovl{(X,d)},x_0)$. For this reason, one usually only considers complete spaces when dealing with Gromov-Hausdorff limits.  In any case, Proposition \ref{GH converge} implies that if a metric space $(X,d)$ supports a $(p,q)$-metric cotype inequality, then the completion $\ovl{X}$ also supports a $(p,q)$-metric cotype inequality. 
\end{remark}

A metric space $(Z,\rho)$ is a \emph{weak tangent} of a metric space $(X,d)$ at a point $p \in X$ if there is a sequence $\{\lambda_k\}$ of positive numbers and a point $p_0 \in Z$ such that the sequence of pointed metric spaces $\{((X,\lambda_k d),p)\}_{k \in \nats}$ pointed Gromov-Hausdorff converges to $((Z,\rho),p_0).$ 

\begin{corollary}\label{weak tan cotype} Suppose that $(X,d)$ supports a $(p,q)$-metric cotype inequality.  If  $(Z,\rho)$ is a weak tangent of $(X,d)$, then $(Z,\rho)$ also supports a $(p,q)$-metric cotype inequality.
\end{corollary}

\begin{proof} By Proposition \ref{scaled bi-Lip invariance}, for any sequence $\{\lambda_k\}$ of positive numbers, the spaces $\{(X,\lambda_k d)\}$ support a $(p,q)$-metric cotype inequality with a uniform constant and scaling function. Proposition \ref{GH converge} now yields the desired result.
\end{proof}

Knowledge of the weak tangents of a space $X$ can now be used to give lower bounds on $q_X$.

\begin{example}\label{cotype 2 ex} Suppose that $X$ is any one of the spaces $\rats$, $\ints$, and $[0,1]$, equipped with the standard metric.  Then $X$ has the Hilbert space $\reals$ as a weak tangent, and of course $X$ is also isometrically embedded in $\reals$.  Since no Banach space has Rademacher cotype less than $2$, and every Hilbert space has Rademacher cotype $2$, Theorem \ref{MN main} and Corollary \ref{weak tan cotype} show that $q_X=2$.  \end{example}

\begin{example}\label{cotype 2 tangent} Suppose that a metric space $(X,d)$ has a weak tangent that contains a rectifiable curve.  Since a rectifiable curve is a bi-Lipschitz image of $[0,1]$, we may conclude from the Example \ref{cotype 2 ex} and Corollary \ref{weak tan cotype} that $q_X \geq 2$. 
\end{example}

\section{Line fitting}

In this section we prove Theorems \ref{1 snowflake} and \ref{p snowflake} using the ideas of the previous sections and the characterization of snowflake spaces, due to Laakso, that is found in the appendix to \cite{TW}.  

\begin{definition}[Laakso]\label{line fitting def} A metric space $(X,d)$ is \emph{line fitting} if for every $c >0$, there is a metric $d_c$ on the disjoint union $X \coprod [0,1]$ with the following properties:
\begin{itemize}
\item there is a number $\lambda_c>0$ such that $d_c|_{X\times X} =\lambda_c d$,
\item $d_c|_{[0,1] \times [0,1]}$ is the standard metric on $[0,1]$, 
\item for each point $t \in [0,1]$, there is a point $x_t \in X$ such that $d_c(x_t,t)<c$.
\end{itemize}
\end{definition}

\begin{theorem}[Laakso]\label{line fitting 1 snowflake} A metric space is a $1$-snowflake if and only if it is line fitting.
\end{theorem}

Given this result, the proof of Theorem \ref{1 snowflake} becomes very similar to the proof of Proposition \ref{GH converge}.

\begin{proof}[Proof of Theorem \ref{1 snowflake}] Throughout, we denote by $|a-b|$ the standard distance between real numbers $a$ and $b$. Suppose that $(X,d)$ is a 1-snowflake.  Let $1\leq p \leq q < 2$, and suppose that $(X,d)$ supports a  $(p,q)$-metric cotype inequality with constant $\Gamma \geq 1$ and scaling function $m_{p,q}(n)$.  We will show that $[0,1]$ also supports such an inequality.  This is a contradiction to Corollary \ref{weak tan cotype}, as discussed in example \ref{cotype 2 ex}. 

Let $n \in \nats$, let $m = m_{p,q}(n) \in 2\nats$, and consider $f \colon \Zmn \to [0,1]$.  Fix $c>0$.  By Theorem \ref{line fitting 1 snowflake}, we may find a metric $d_c$ on $X \coprod [0,1]$ as in Definition \ref{line fitting def}.  Hence there is a function $\til{f} \colon \Zmn \to X$ with the property that for any $\ep \in \Zmn$,
$$d_c(\til{f}(\ep),f(\ep)) < c.$$ 
From this, the triangle inequality, and the other properties of $d_c$, we see that 
$$\E_{\ep \in \Zmn} \sum_{j=1}^n |f(\ep)-f_j(\ep)|^p  \leq n(4c)^p +2^{p}\lambda_c^p\E_{\ep \in \Zmn} \sum_{j=1}^n d(\til{f}(\ep), \til{f_j}(\ep))^p.$$
A similar argument shows that 
$$\E_{\ep \in \Zmn} \E_{\del \in \{-1,0,1\}^n} d(\til{f}(\ep),\til{f}_{\del}(\ep))^p \leq \lambda_c^{-p}(4c)^p + 2^{p}\lambda_c^{-p}\E_{\ep \in \Zmn} \E_{\del \in \{-1,0,1\}^n} |f(\ep)-f_\del(\ep)|^p.$$
Applying the $(p,q)$-metric cotype inequality on $X$ to the function $\til{f}$ and combining with the previous estimates produces 
\begin{align*}\E_{\ep \in \Zmn} \sum_{j=1}^n |f(\ep)-f_j(\ep)|^p \leq (4c)^p ( n + & 2^{p}\Gamma^pm^pn^{1-(p/q)})  \\ & + 4^{p}\Gamma^pm^pn^{1-(p/q)} \E_{\ep \in \Zmn} \E_{\del \in \{-1,0,1\}^n} |f(\ep)-f_\del(\ep)|^p.\end{align*}
Letting $c$ tend to $0$ completes the proof.
\end{proof}

\begin{proof}[Proof of Theorem \ref{p snowflake}] We first claim that as $(X,d)$ is an $s$-snowflake, there is an $(1/s,L)$-snowflaking embedding, $L \geq 1$, from a $1$-snowflake space into $(X,d)$.  The desired result then follows from Theorem \ref{1 snowflake} and Proposition \ref{sf invariance}.

By definition, the space $(X,d)$ is $L$-bi-Lipschitz equivalent, $L \geq 1$, to a metric space $(Y,d_Y)$ with the property that $(Y,d_Y^s)$ is also a metric space.  This immediately implies that there is a $(1/s, L)$-snowflaking embedding of $(Y,d_Y^s)$ into $(X,d)$. If $(Y,d_Y^s)$ is not a $1$-snowflake, then it is bi-Lipschitz equivalent to a metric space $(Z,d_Z)$ such that $\rho=d_Z^t$ is also a metric on $Z$ for some $t>1$.  Since $st \geq 1$, the distance $\rho^{1/(st)}$ is again a metric on $Z$. Moreover, $(Z,\rho^{1/(st)})$ is an $L^{st}$-metric space that is bi-Lipschitz equivalent to $(X,d)$. This is a contradiction with the assumption that $s=s_X$, as $st>s$. 
\end{proof}

\bibliographystyle{plain}
\bibliography{MetricCotypePost}
\end{document}